\def\cyrb{\bf}
\def\cit{\it}
\def\O{\Omega}
\def\g{\gamma}
\def\a{\alpha}

\def\g{\gamma}
\def\e{\varepsilon}
\def\f{\varphi}
\def\l{\lambda}
\font\Bbb=msbm10
\def\R{\hbox{\Bbb R}}

\def\Z{\hbox{\Bbb Z}}
\def\N{\hbox{\Bbb N}}

\def\An{\mathop{\rm An}} 
\def\Err{\mathop{\rm Err}}  
\def\const{\mathop{\rm const}}  
\def\tg{\mathop{\rm tg}}   
\def\x{{\hbox{\bf x}}}
\def\y{{\hbox{\bf y}}}  
\def\z{{\hbox{\bf z}}}  
\tolerance 10000
                                    
{\rightline{\it F. Petrov}}  

{\bf On a number of rational points on a convex curve}
\smallskip
{
Let $\g$ be a bounded
convex curve on a plane. Then
$\sharp (\g\cap (\Z/n)^2)=o(n^{2/3})$.
It streghtens the classical result of Jarn\'\i k
[J] (an upper estimate $O(n^{2/3})$) and 
disproves a conjecture of Vershik
on existence of the so-called
{\it universal Jarn\'\i k curve}.
} 
\smallskip

\footnote{}{MSC 11H06, {\it Key words:}: 
convex curve, lattice points, 
affine length}
\footnote{}{
The research is supported by grants
NSh.2251.2003.1 (Nauchnaya shkola) 
and RUM1-2622-ST-04 (CRDF).}

In his famous paper [J] Jarn\'\i k proved
(among other reuslts), that the maximal
possible number of integer points, which may lie
on a strictly convex plane curve of length $N$ grows
like $cN^{2/3}$ (the exact constant $c$ was also
computed in [J] and equals $3(2\pi)^{-1/3}$). 
In other words, the number of nodes of a lattice
$L_N:=(\Z/N)^2$, which lie on a strictly convex
curve $\g$ of length 1, does not exceed $cN^{2/3}$;
and for any $N$ there exists such strictly convex
curve $\g^{(N)}$ of length 1, that
$$k(\g^{(N)},N):=\sharp(\g^{(N)}\cap L_N)\ge cN^{2/3}\eqno (*)$$

So, the natural question arises: does there exist a
{\it universal} curve $\g$, for which (*) 
holds for infinite number of positive integers
$N$? This question is formulated
by A. M. Vershik;
in the paper
[P] it was attributed to J.-M. Deshoullirs 
and G. Grekos. The                   
conjectire about the existence of a universal curve
was formulated in [P]. 
Indeed, the methods of
papers [Ve],[Ba] (see below) show 
the nature of the typical convex
lattice polygon, which, 
on the first glance,
supports such conjecture.
However, this conjecture fails and here we disprove it.

H.~P.~F.~Swinnerton-Dyer [SD] has proved an estimate
$k(\g,N)\le cN^{3/5+\e}$ for
$\g\in C^3$, and E.~Bombieri and J.~Pila [BP] 
proved an estimate
$k(\g,N)\le cN^{1/2+\e}$ ($\e>0$ is arbitrary)
for infinitely smooth $\g$.
Here we do not mention further results in this
direction in terms of smoothness, curvature
conditions and other restrictions on a curve.

A. M. Vershik and I. Barany ([Ve], [Ba]) 
investigated limit shapes of large
random polygons with vertices on a shallow lattice.
The answers bring out the connection
with affine geometry. Namely, let
$l_a(\g)$ denote an affine length of a curve $\g$ 
(an integral of cubic root of curvature by natural
parameter). It appears that the number
of polygons with vertices in nodes of
$L_n$, lying in a small
neighborhood of given curve $\g$,
grows like $e^{c\cdot l_a(\g) n^{2/3}}$, 
and the number of their vertices --- like $c\cdot l_a(\g)n^{2/3}$
(remarkably, it holds both for
maximal and typical number of vertices,
only constants differ). 
Polygons with vertices in nodes of $L_n$, 
which lie inside a given convex polygon,
concentrate near the closed convex curve
with maximal possible affine length ([Ba]). 
This curve is nothing else
but the union of some parabola arcs,
inscribed in polygon angles. Hence the number
of $L_N$ nodes on such curve 
does not exceed $C\cdot N^{1/2}$. So, the typical
curve is not universal. Here we prove that
in fact the universal curve does not exist:
for any bounded strictly convex curve an estimate
$k(\g,n)=o(n^{2/3})$ holds. 

I am deeply grateful to A. M. Vershik for posing the
problem, numerous helpful discussions, attention and
manuscript editing. I am also grateful to A. Gorodnik
for valuable consultations on equidistribution and
to S. Duzhin, who pointed me out the book [F].

\smallskip

{{\bf \S 1. Definitions and denotions.}}
\smallskip
Fix a cartesian coordiante system on a plane.

Let $S(F)$ denote the doubled area of a polygon 
$F$; $\x\times\y$ denote the pseudoscalar product
of vectors $\x$ and $\y$
(i.e. an oriented area of a parallelogram,
based on these vectors).

Fix a triangle $ABC$, oriented in such manner that
$S=S(ABC)=+\overline{AC}\times \overline{CB}$.

Define following concepts:

1. $\An=\An(ABC)$ is an angle with
vertex in the origin and sides collinear with rays
$AC$ and $CB$ (value of this angle equals
$\pi-\angle ACB$). We consider angle as a set of {\it vectors}
$\x\in \R^2$, which, being drawn from the origin,
lie inside $\An$. In other words, $\An$
is a set of vectors 
$\{t_1\cdot\overline{AC}+t_2\cdot\overline{CB}|\, t_1,t_2\ge 0\}$

2. Given a vector $\x$, define its
{\it girth} as $$[\x]=(\x\times \overline{CB}+\overline{AC}\times\x)/S.$$
Girth is a linear function
of a vector
$\x\in \R^2.\,$
Note that $[\overline{AC}]=[\overline{CB}]=1,\,[\overline{AB}]=2$.
We also define a girth of a segment 
$PQ$ as a value $[PQ]=|[\overline{PQ}]|$.

3. Define $ABC-$radius of arbitrary triangle
as a product of its sides girths, divided by
quadruple area (it appears to be
the usual circumradius if to take sides lengths
instead girths). Note that
$ABC$-radius of triangle
$ABC$ equals $S^{-1}$.

4. Let $\g:AC_1C_2\dots C_kB$ be a strictly
convex broken line. We call it 
$(AB,C)$-broken line, if all its
vertices lie inside tiangle
$ABC$. If, additionally, all 
intermediate vertices 
$C_i (i=1,\,2,\,\dots,k)$ lie in nodes of a lattice
$L_n=({1\over n}\Z)^2$, we call $\g$ an
$(AB,C;n)$-broken line.

5. Let $(AB,C)$-broken line 
$\g=AC_1C_2\dots C_kB$ be inscribed in an 
$(AB,C)$-broken line $\g_1=AD_1D_2\dots D_{k+1}B$
(i.e. points $C_i$ lie on respective segments
$D_iD_{i+1}\,(i=1,\,2,\,\dots,\,k)$).
Define a {\it generalized affine length of a broken line
$\g$ with respect to $\g_1$} as
a value
$$
l_A(\g:\g_1):=\sum_{i=0}^{k} 
S(C_iD_{i+1}C_{i+1})^{1/3}\, (C_0=A,\,C_{k+1}=B),
$$
and a
{\it generalized affine length of a broken line $\g$} as a value
$$
l_A(\g)=\sup_{\g_1} l_A(\g:\g_1),
$$
where the supremum is taken by 
all $(AB,C)$-broken lines $\g_1$,
circumscribed around $\g$.

The generalized affine length of $(AB,C)$-broken line
$\g$ may be also defined (up to some factor)
as a supremum of affine lengthes of smooth curves,
circumsribed around $\g$ and contained inside triangle
$ABC$. In classical affine geometry, the affine length of a broken
line equals 0, so we have to add justifying "{\it generalized}"
in our definition.

\smallskip
{\bf \S 2. Preliminary statements.}
\smallskip
At first, we need an upper estimate
for a numebr $k$ of intermediate vertices of an
$(AB,C;n)$-broken line for large $n$.
We use the known fact that an area of convex
$k$-gon with integer vertices is not less than
$(8\pi^2)^{-1}k^3\ge (k/5)^3\,(k\ge 3)$.
From this,
$$
k\le \max(3,5(Sn^2)^{1/3})\,  \eqno(1)
$$

So the maximal quantity
of nodes of lattice 
$L_n$, which lie in a triangle of area
$S$, grows not faster than $c(Sn^2)^{1/3}$.

At fact, it grows exactly as 
$c(Sn^2)^{1/3}$. 
An example of $(AB,C;n)$-broken line,
whihc has about $c(Sn^2)^{1/3}$ vertices may be
constructed as follows:
arrange vectors of a set 
$\An\cap L_n$ in girth increasing order.
Take $c(Sn^2)^{1/3}$ vectors with minimal
girths and construct a convex broken line,
for which these vectors are vectors of 
edges. For small enough $c$ (say, $c=1/100$)
this broken line may be shifted inside
triangle $ABC$ in such way that its
vertices will still lie in $L_n$ nodes 
and will form, together with vertices 
$A$ and $B$, an $(AB,C)$-broken line. 
This construction is quite analogous 
to Jarn\'\i k's construction of broken line
with given length and maximal number
of integer points, and may be consider as its
affinne generalization.
This example is needed not before paragraph 4.
 
The idea of our approach to the statement
mentioned in abstract is as follows.
Draw tangents to our curve in points, which lie
in nodes of $L_n$. For large $n$,
we get many small triangles,
and their union contains our curve.
It appears that if our curve contains quite many
points on these lattices, then
the sum of cubic roots of areas of these
triangles decreases, with
any new $n$, by some factor, which is impossible.
For proving lemma 4 we need some technical statements
form this paragraph.

%


Now consider all vectors of a set
$\Z^2\cap \An$ and arrange them by girths
(in ascending order).
Consider $k$ such vectors $\z_1,\,\z_2,\,\dots,\z_{k}$
with the least girths. 

An array of vectors
$\z_1,\,\z_2,\,\dots,\,\z_{k}$ has the following structure:
it contains all vectors, which lie
in triangle $OPQ=\An\cap\{\z:[\z]<r=[\z_{k}]\}$,
and some vectors, which girth equals $r$.
We have
$$
k={1\over 2}S(OPQ)+o(r^2)={S\over 2} r^2+o(r^2) 
$$
Further, 

$$
\sum [\z_i]=\int_{OPQ}[\x]d\x+o(r^3)
=r^{-1}S^{-1}\int_{OPQ}(S(OPX)+S(OQX))d\x+o(r^3)=
{S\over 3}r^3+o(r^3).
$$

Hence 

$$
\sum [\z_i]\ge c S^{-1/2} k^{3/2} +o(k^{3/2}),\,c=2\sqrt{2}/3 \eqno(2)
$$

The following elementary lemma
is a basepoint of further considerations.

{\bf Lemma 1.} Let points
$P,\,R$ be chosen on the sides $AC$ and
$BC$ of a triangle $ABC$ respectively, and a point
$Q$ --- on a segment $PR$.
Then:

$1^{\circ}$. $S(AQP)^{1/3}+S(BQR)^{1/3}\le S^{1/3}.$
Moreover, there exists a function 
$\e_1(\e)$, which tends to zero with 
$\e\rightarrow 0$, such that if  \smallskip
$\Err:=1-(S(APQ)/S)^{1/3}-(S(BQR)/S)^{1/3}<\e [AQ]$,
then

{$2^{\circ}$.} $[AP]:[PQ]\in (1-\e_1,1+\e_1)$  

{$3^{\circ}$.} $r(AQP)\in S^{-1}(1-\e_1,1+\e_1)$, 
where $r(APQ)$ denotes $ABC$-radius of a triangle $APQ$.

\item{$4^{\circ}$.} For any vector
$\x\in \An(APQ)$ we have $[\x]_{APQ}:[\x]\in (AP:AC)\cdot(1-\e_1,1+\e_1)$
(here $[\x]_{APQ}$ denotes an
$APQ$-girth of a vector $\x$).

{\bf Proof.} 
We have 

$
\Err=1-(S(APQ)/S)^{1/3}-(S(BRP)/S)^{1/3}= 
1-({AP\over AC}\cdot {PQ\over RP}\cdot {RC\over CB})^{1/3}-
({PC\over AC}\cdot {QR\over RP}\cdot {BR\over BC})^{1/3}=
$

$=\left(
({AP\over AC}+{PQ\over RP}+{RC\over CB})/3- 
({AP\over AC}\cdot {PQ\over RP}\cdot {RC\over CB})^{1/3}
\right)+
\left(
({PC\over AC}+{QR\over RP}+{BR\over BC})/3-
({PC\over AC}\cdot {QR\over RP}\cdot {BR\over BC})^{1/3}
\right)\ \ (6)
$

Two last brackets have a form
$(x+y+z)/3-(xyz)^{1/3}$
(and by AM-GM inequality, it implies p.1 1 of lemma).
We have
$$2(x+y+z)-6(xyz)^{1/3}=(x^{1/3}+y^{1/3}+z^{1/3})\left((x^{1/3}-y^{1/3})^2+
(y^{1/3}-z^{1/3})^2+(z^{1/3}-x^{1/3})^2\right).$$
So, if $\max(x,y,z)\ge (1+\delta) \min(x,y,z)$,
then $(x+y+z)/3-(xyz)^{1/3}\ge c(\delta) (x+y+z)$,
where $c(\delta)$ is a positive function, which tends
to 0 with$\delta\rightarrow 0$.

Applying this observation
to the first summand in RHS of $(6)$ 
($x=AP/AC,\,y=PQ/RP,\,z=RC/CB$), 
we get the following implication: if $\Err\le \e [AQ]$,
then (since

$[AQ]=[AP]+[PQ]=AP/AC+(PQ/PR)\cdot[PR]\le (AP/AC)+2(PQ/PR)$),
the estimate $\Err\le 2\e(x+y+z)$ holds, and therefore
$\max(x,y,z):\min(x,y,z)<1+\e_1,$
where $\e_1$ is chosen such that
$c(\e_1)>2\e$. Such $\e_1(\e)$ may be chosen, and
it tends to zero with $\e$.

Further,
$$
[PA]=x,\,[PQ]=(PQ/PR)\cdot [PR]=y([PC]+[CR])
=y(1-x+z)=[PA]\cdot(y/x)\cdot(1-x+z)
$$
We have 
$$
y/x\in (1-\e_1,1+\e_1),\, 
1-x+z=1+z(1-x/z)\in (1-z\e_1,1+z\e_1)\subset (1-\e_1,1+\e_1),
$$
hence
$$
[PA]:[PQ]\in((1-\e_1)^2,(1+\e_1)^2)\subset (1-\e_2,1+\e_2),\,\e_2:=2\e_1+\e_1^2,
$$
which proves p.2
(with $\e_2$ instead $\e_1$).

For proving p.3 we note that
$$
S\cdot r(ABC)=[AP]\cdot[PQ]\cdot({[AP]+[PQ]\over 2})\cdot (S/S(APQ))=
{[AP]\cdot[PQ]\cdot([AP]+[PQ])\over 2xyz}\in (1-\e_3,1+\e_3) 
$$
where $\e_3(\e)\rightarrow 0$ if $\e\rightarrow 0$
(since mutual quotients of numbers
$[AP]=x,\,[PQ],\,y,\,z$ lie in a segment $(1-\e_2,1+\e_2)$).

We pass to the p.4.
Let $\x=a \overline {AP}+b\overline {PQ},\,\x\in \An(APQ)\,(a,\,b\ge 0)$
Then $[\x]_{APQ}=a+b,\,[\x]=a[AP]+b[PQ]=(a+b)[AP]+b([PQ]-[AP])=
(AP/AC)(a+b)(1+{b\over a+b}\cdot([PQ]/[AP]-1))$.
The statemnet of p.4 follows.
\smallskip
{\bf Note.} Point 1 of lemma may be founded, for instance,
in book [F] (p. 391).  
\smallskip
{\bf Corollary.} Applying the statemnt of p.1.
many times, we get the following important fact:
generalized affine length of any $(AB,C)$-broken line $\g$
satisfies an inequality $l_A(\g)\le S^{1/3}$.
Moreover, generalized affine length does not increase
with adding new vertices to a broken line.
\smallskip
Now we formulate a statement on asymptotic
distribution of integer points on a surface
$ab-cd=\const$.

{\bf Lemma 2.} Consider pairs of vectors 
$(\x,\,\y): \x,\,\y\in \An\cap \Z^2$, for which
$\x\times \y=m\ne 0$ ($m\in \Z$ is a constant).
For each such pair, consider a {\it special} point
$([x],[y])\in [0,\infty)^2.$ Then special
points are equidistributed in a first quadrante in a following sense:
for any bounded domain $\O\subset (0,\infty)^2$
with piecewise-smooth boundary, the number of special
points in a domain $N\Omega$ has an asymptotics
(by $N\rightarrow \infty$)
$$
c(m)S\cdot N^2S(\Omega)+o(N^2),
$$
where $c(m)$ is a constant, which depends on $m$
(namely, $c(m)=(2\zeta(2))^{-1}\sigma(m)/m$,
$\sigma(m)$ is a sum of positive integer divisors of 
$m$).
\smallskip

The proof of lemma 2 is given in appendix. Alex Gorodnik 
pointed out that these statements follow from known general
results (for instance, [EM]).

\smallskip
{\bf Corollary.} Consider triangles
$PQR$ such that 

$1^{\circ}.\,\overline{PQ},\,\overline{QR}\in \An\cap \Z^2$

$2^{\circ}.\,[PR]\le M(n/S)^{1/3}$ 

$3^{\circ}.\,r(PQR)\in 2nS^{-1}(t_1,t_2)\,(0<t_1<t_2)$

$4^{\circ}.\, S(PQR)\le m\,(m\in \N)$.

The number $N(m,M,t_1,t_2)$ of such triangles
may be bounded by above as follows ($n\rightarrow \infty$):

$$
N(m,M,t_1,t_2)\le c(m)(t_1-t_2)^{3}(nS^2)^{1/3}+ o(n^{2/3}).
$$
\smallskip
{\bf Proof.} Just apply lemma 2 for domains
$\Omega(m',M,t_1,t_2)=\{(p,q):0<p,q<M,\, pq(p+q)\in 4m'(t_1,t_2)\}$
and $N=(n/S)^{1/3}$ ($m'=1,\,2,\,\dots,\,m$) and use a fact
that domain $\Omega(m',\infty,0,1)$ has a finite area.
\smallskip
{\bf Lemma 3.} 
Consider convex quadrilateral
$PSTR$ such that 
$\overline{PS},\,\overline{ST},\,\overline{TR}\in \An$.
Let $Q$ be a point on its side $ST$.
Then the value of
$r(PQR)$ lies between numbers
${[PQ]\over [PS]}r(PQS)$
and ${[QR]\over [RT]}r(QTR)$.
\smallskip


{\bf Proof.} We have
$S(PQR)={[QR]\over [SQ]}S(PQS)+
{[PQ]\over [QT]}S(RQT)$. This equality may \smallskip
seem unexpected: RHS a priori depends on a choice of linear
function $[\x]$. It may be proved,
for instance, as follows: when points $S$ and $T$ move on fixed rays
$QS$ and $QT$, nor left, neither right sides of
equality do not change. So, we may assume that
$[SQ]=[QR]$ and $[QT]=[PQ]$. 
First equality means that the line, which joins
$Q$ and a midpoint of $RS$, is parallel to a vector
$\overline{AC}+\overline{BC}$.
Analogously, second equality means that the line,
which joins $Q$ and the midpoint of $PT$ is
parallel to the same vector.
So, the point $Q$ lies on a Gauss line of a quadrilateral
$PRTS$, and therefore satisfies an equation of this line
$S(PXR)+S(TXS)=S(SXP)+S(RXT)$
(of course, to get the Gauss line
we must consider oriented areas in this equation), 
q.e.d.
From here
$$
2r(PQR)=
{
[PQ]\cdot [QR]\cdot ([PQ]+ [QR])\over
{[QR]\over [SQ]}S(PQS)+
{[PQ]\over [QT]}S(RQT)
}=
{[PQ]+ [QR] 
\over
{[PS]\over 2r(PQS)}
+
{[TR]\over 2r(QTR)} 
}
$$
As it is known, the quotient ${x+y\over x'+y'}$ lies
between ${x\over x'}$ and ${y\over y'}$. 
Applying this for
 $x=[PQ],\,y=[QR]$
and $x'={[PS]\over r(PQS)},
\,y'={[TR]\over r(QTR)}$,
we get desired statement.

\smallskip

{\bf \S 3. Main part}
\smallskip

{\bf Theorem.} Let $\g$ be a bounded convex
curve. Denote $k(\g,n)=\sharp(\g\cap L_n)$. Then 
$k(\g,n)=o(n^{2/3})$. 
\smallskip
{\bf Lemma 4.} For any $c>0$ there exists such
a number $a(c)>0$, that for any
triangle $\triangle ABC$ for large enough
$n>N(c,\triangle ABC)$ the following statement 
holds: any $(AB,C;\,n)$-broken line $\g$,
which have at least 
$\ge c(n^2S(ABC))^{1/3}$ vertices, have not too much
generalzied affine length: $l_A(\g)\le (1-a(c))S(ABC)^{1/3}$.

At first, we show how the lemma 4 implies a theorem.

Asuume that the theorem is not
valid and for some bounded convex
curve $\g$ we have
$k(\g,q_n)\ge cq_n^{2/3}$ for some increasing
sequence of positive integers $q_1<q_2<q_3\dots$.
Without loss of generality, curve
$\g$ joins points $A$ and $B$, 
and lies inside triangle $ABC$
(every bounded convex
curve may be partitioned onto
finite number of such curves
). Let $S(ABC)=1$. Define 
$\g_n$ as an $(AB,C)$-broken line, inscribed in $\g$,
for which the set 
$\g\cap (\cup_{i=1}^n L_{q_i})$ is a set of intermediate
vertices.

We fix a support line in each point of
curve $\g$. Draw these lines in vertices of
$\g_n$, we get an
$(AB,C)$-broken line $\g_n'$, circumscribed
around $\g_n$. Denote by
$\triangle_1,\,\triangle_2,\,\dots,\triangle_k$,
triangles formed by intersecting support
lines in neighbour vertices of $\g_n$.
The line $\g$ lies inside the union
$\cup_{i=1}^k \triangle_i$. Denote $S_i=S(\triangle_i)$.
Let $q=q_m$ be such large integer, that
for any triangle $\triangle_i\,(1=1,\,2,\,\dots,k)$
the alternative of lemma 4 holds:
either

\item{(1)} 
$\sharp(\g\cap \triangle_i\cap L_q)\le (c/2)\cdot S_i^{1/3}\cdot q^{2/3}$, 
or

\item{(2)} $l_A(\g_m\cap\triangle_i)\le (1-a)S_i^{1/3}$.

Note once more, that number $a$ depends only
on $c$. Denote by $M_1$ the set of indecies $i$,
for which the case (1) holds, and by $M_2$ --- the set of ther
indecies $i$ (for them, case (2) holds).
Note that
$$
c(\sum_{i=1}^k S_i^{1/3}) q^{2/3}\le 
cS(ABC)^{1/3}q^{2/3}=cq^{2/3}\le \sharp(\gamma\cap L_{q})\le
\sum_{i=1}^k \sharp(\gamma\cap L_q\cap \triangle_i)
=\sum_{i\in M_1}
\sharp(\gamma\cap L_q\cap \triangle_i)+
$$

$$
+\sum_{i\in M_2}
\sharp(\gamma\cap L_q\cap \triangle_i)
\le(c/2)(\sum_{i\in M_1} S_i^{1/3}) q^{2/3}+
5(\sum_{i\in M_2} S_i^{1/3}) q^{2/3}\le 
(c/2)(\sum_{i=1}^k S_i^{1/3})q^{2/3}+
5(\sum_{i\in M_2} S_i^{1/3}) q^{2/3}, 
$$
hence
$$\sum_{i\in M_2} S_i^{1/3}\ge {c\over 10}(\sum_{i=1}^k S_i^{1/3}).$$
Furthermore, 
$$l_A(\g_{m}:\g_m')\le \sum_{i=1}^k 
l_A(\g_m\cap \triangle_i)=
\sum_{i\in M_1}+\sum_{i\in M_2}\le
$$

$$
\le
\sum_{i\in M_1} S_i^{1/3}+(1-a)\sum_{i\in M_2} S_i^{1/3}
\le (1-{ac\over 10})(\sum_{i=1}^k S_i^{1/3})=(1-{ac\over 10})l_A(\g_n:\g_n').
$$
From here we may deduce that $\lim_{n\rightarrow \infty} l_A(\g_n:\g_n')=0$,
which contradicts to our
assumptions
(from which $l_A(\g_n:\g_n')\ge c/5$).

So, it suffices to prove lemma 4.

Let $\g=C_0C_1C_2\dots C_k C_{k+1}\,(C_0=A,\,C_{k+1}=B,\,
C_i\in L_n (i=1,\,2,\,\dots,\,k))$ be an $(AB,C;n)-$broken line, and
$k\ge c S^{1/3}n^{2/3}$. Let's
fix supposrt lines $l_i$
in points $C_i$ ($i=1,\,2,\,\dots,\,k$),
define also lines
$l_0=AC$ and $l_{k+1}=BC$. 
Put $l_i\cap AC=B_i\,(i=1,\,2,\,\dots,\,k+1)$,
$l_i\cap l_{i+1}=D_i\,(i=0,\,1,\,\dots,\,k)$.

\medskip

For optimal (in sense of generalized affine length definition)
choice of support lines we have
$$
l_A(\g)=\sup \sum_{i=0}^k S(C_iD_iC_{i+1})^{1/3},
$$
We have 
$$
S(ABC)^{1/3}-l_A(\g)=\sup \sum_{i=1}^k x_i,\, 
x_i=S(AC_{i+1}B_{i+1})^{1/3}-S(C_iC_{i+1}D_i)^{1/3}-S(AC_iB_i)^{1/3}.
$$
Numbers $x_i$ are non-negative by p.1 of lemma 1.
Our goal is to find 
$N\ge c_1n^{2/3}S(ABC)^{1/3}$ indecies $i$, for which
$x_i\ge c_2 [C_iC_{i+1}] S(ABC)^{1/3}$ 
(with some constants $c_1$ and $c_2$, which
depend only on $c$). 
If we succeed, the lemma 4 will be proved because
of estimate
(2) for the sum of $N$ least girths.

Note also that if we find
$N_1\ge c_3n^{2/3}S(ABC)^{1/3}$ indecies
$i$, for which $x_i\ge c_4 [C_iC_{i+1}] S(AC_{i+1}B_{i+1})^{1/3}$,
it is also enough for our goal. Indeed, using inequality
$(5)$, we have an estimate like

$100 n^2S(AC_{i+1}B_{i+1})\ge i^3$, so
for $i\ge N_1/2$ we have $S(AC_{i+1}B_{i+1})^{1/3}\ge c_6S(ABC)^{1/3}$,
and hence for at least $N_1/2$ indecies $i$
desired estimate $x_i\ge c_7 [C_iC_{i+1}] S(ABC)^{1/3}$
holds.

Applying lemma 1 at first stage for triangle
$ABC$ and points $P=B_{i+1},\,Q=C_{i+1}$,
and then for triangle $AC_{i+1}B_{i+1},\,P=B_i,\,Q=C_i$,
we see, that it suffices for some
$\e_0(c)$ to find $N_1\ge c_3n^{2/3}S(ABC)^{1/3}$
indecies $i$, for which either
$\max([C_iD_i],[D_iC_{i+1}]):\min([C_iD_i],[D_iC_{i+1}])>1+\e_0$
or $r(C_iD_iC_{i+1})\notin S^{-1}(1-\e_0,1+\e_0)$.
We call the index $i$, satisfying at least one 
of these two conditions, $\e_0$-nice. 
Lemma 3 (for $PSTR=C_iD_iD_{i+1}C_{i+2},\, Q=C_{i+1}$
implies the following

{\bf Statement.} Assume that indecies
$i$ and $i+1$ are not
$\e_0$-nice. Then for some
$\e_1(\e_0)$ we have
$r(C_iC_{i+1}C_{i+2})\in 2S^{-1}(1-\e_1,1+\e_1)$,
and $\e_1$ tends to 0 together with $\e_0$.

So, it suffices to prove that for
$\e>0$ the number $N_{\e}$ of indecies $i$,
for which $r(C_{i}C_{i+1}C_{i+2})\in 2S^{-1}(1-\e,1+\e)$,
admits an upper bound $N_{\e}\le \e_1 (Sn^2)^{1/3}$,
where $\e_1$ tends to 0 together with $\e$.

Note that $\sum S(C_iC_{i+1}C_{i+2})^{1/3}\le 2S^{1/3}$
(sums by odd and even $i$ do not exceed $S^{1/3}$
by corollary of lemma 1), so the number of indecies
$i,$ for which
$S(C_iC_{i+1}C_{i+2})\ge mn^{-2}$,
does not exceed (by Chebyshev inequality) 
$2m^{-1/3} (SN^2)^{1/3}$.
Moreover, since for the sum of girths we have
$\sum_{i=0}^k [C_iC_{i+1}]\le 2$,
we may aslo assume that
$[C_iC_{i+2}]\le M(n^2S)^{-1/3}$
(Chebyshev inequality again) for some large $M$.

Now it suffices to use the corollary
of lemma 2.
\smallskip
{\bf Note.} The statement of a theorem
may be slightly strengthened by letting
$L_n=({1\over n}\Z)^2+\x_n$,
where $\x_n$ is some vector of a shift,
and chosing $n$ not necessary integral.
\smallskip
{\cyrb \S 4. About possible
number of integer points on a convex curve}
\smallskip

The following example was constructed in 
[P]: there exists a strictly
convex bounded curve $\g$, for which
$k(\g,q_n)\ge c_n q_n^{2/3}$, where
$q_1<q_2<\dots$ is arbitrary fast increasing 
integer sequence, and coefficients $c_n$
decrease as $K^{-n} (K>1$ is some explicit constant).
Here we stengthen this result: it suffices to suppose that
$\sum c_n<\infty$. But our method is quite rigorous,
and it seems probable that this condition 
may be reduced to some weaker one
(may be even to the necessary condition
$\lim c_n=0$).

{\bf Theorem 2.} Let $\sum c_k<\infty$ be a convergent positive
series; $M\subset \N$ be an infinite subset of
positive integers. Then there exists a sequence
$q_1<q_2<\dots,\,q_i\in M$ of positive integers
and a strictly convex bounded curve $\g$ such that
$k(\g,q_n)\ge c_n q_n^{2/3}$

{\bf Proof.} Consider a semicircle of length 
$\sum c_i$. Partion it onto arcs of lengthes $c_i$,
denote by $A_{i}$ and $A_{i+1}$ endpoints of arc of length 
$c_i$. Draw tangents to a semicircle in endpoints
of these arcs. Each arc lies in some triangle
$A_iB_iA_{i+1}$ ($B_i$ is a point,
in which tangents in $A_i$ and $A_{i+1}$ intersect).  
An area of triangle $A_iB_iA_{i+1}$ is not less than
$Cc_i^3$, where  $C$  does not depend on $i$. 
Without loss of generality, $C>100$.
(else make a homothety, which increase areas
in $100/C$ times). Now take a single triangle
$A_iB_iA_{i+1}$ and construct such an
$(A_iA_{i+1},B;q_i)$-broken line, that number
of intermediate vertices of this broken line
is not less than $c_iq_i^{2/3}$
(it is possible for large enough
$q_i$, as it was discussed in paragraph 1).    
The union of all these broken lines is desired
curve $\gamma$ (to be more precise, we must
replace straigth edges of broken lines to 
some smoothe strictly convex curves for making $\g$ 
strictly convex).  
 \medskip

\centerline{\bf Appendix. Distribution of integer points on a surface 
$ab-cd=\const$.}
\smallskip      

Consider bounded domains
$\O_1,\,\O_2\subset \R^2$. 
with piecewise-smooth boundary.
We study the asymptotics of a quantity
$M(\O_1,\O_2;\,n)$ of pairs of vectors
${\x}_1\in n\O_1\cap \Z^2,\,
{\x}_2\in n\O_1\cap \Z^2$ such that 
${\x}_1\times {\x}_2=1$ (here ${\x}_1\times {\x}_2$ 
denotes an oriented area of parallelogramm
based on vectors ${\x}_1,\,{\x}_2$).
At first, we consider the case of triangles
$$\O_i=\{(x,y):0<y<x\le a_i\}\,(i=1,2,\,a_i>0)$$
In this case, the problem reduces to the
following question: how many solutions does 
an equality
$$x_1y_2-y_1x_2=1\eqno(a1)$$ 
with conditions
$$0<y_1<x_1<na_1,\,0<y_2<x_2<na_2\eqno(a2)$$
have?
It is well-known, that for fixed coprime
$x_1,\,x_2$ an equation $(*)$ has a unique
solution in $y_1,\,y_2$,
for which $0<y_1\le x_1,\,0<y_2\le x_2$. 
Cases of equality $y_1=x_1$ or $y_2=x_2$ 
may be realized only if
$x_1=1$ or $x_2=1$, i.e. for at most
$C(a_1,a_2)\cdot n$ variants. 
The number of integer points
$(x_1,\,x_2)$ with coprime coordinates in the rectangle
$0<x_1<na_1,\,0<x_2<na_2$ is
$$\zeta(2)^{-1}n^2a_1a_2+o(n^2),$$
hence for the case (*) we have asymptotics
$$
M(\O_1,\O_2;\,n)=\zeta(2)^{-1}n^2a_1a_2+o(n^2). 
$$
In next, we need the following reformulation of gotten result.
Denote by $l(\O,\f)$ 1-dimensional Lebesgue
measure of the intersection of domain 
$\O$ and a line $y=\tg \f\cdot x$
(the line which form angle $\f$ with X-axis).
Then
$$
M(\O_1,\,\O_2;\,n)=\zeta(2)^{-1}\int_0^{\pi} 
l(\O_1,\f)\cdot l(\O_2,\f) d\f\    n^2+o(n^2)
\eqno(a3)
$$

A triangle $OAB$ with area $1/2$ and integer vertices $A$
and $B$ is called a {\it basic} triangle.

Note that formula (a3)
holds also for domains,
which may be gotten from (*) by
$Sl(2,\Z)$-element affine action.
In other words, formula (a3) holds for triangles,
which may be gotten from some basic triangle
by homotheties.

Our further plan is approximation
of quite generic domains by
unions of such "basic" domains,
almsot disjoint in central projection
to unit circle.

Fix a number $\e>0$. The basic triangle
$OAB$ is called {\it $\e$-suitable},
if

\item{$1^{\circ}.$} $|OA/OB-1|<\e$

\item{$2^{\circ}.$} $\angle AOB<\e$.

We need the following 

{\bf Lemma.} Almost every (in the sense of
Lebesgue measure on a unit circle) 
ray, arising from the origin,
intersects interiors of infinitely many 
$\e$-suitable basic triangles.

{\bf Proof.} Without loss of generality,
the ray is defined as $0<y=\a x,\,0<\a<1$.
Consider a continued fraction for 
$\a$: $\a={1\over a_1+{1\over a_2+\dots}}$.
For almost all $\a$ elements
$a_i$ are unbounded (it follows from Gauss-Kuzmin
formula, but may be gotten easier, see [Kh]).
In terms of convergents
${p_k\over q_k} (k=1,\,2,\,\dots)$ 
it means that the ratio $q_{k+1}/q_k$ of neighbour
denominators is unbounded
$q_{k+1}/q_k$
(since $q_{k+1}=a_kq_k+q_{k-1}$).
Every pair of neighbour convergents
corresponds to some basic triangle
$OAB$ ($A=(q_k,p_k),\,B=(q_{k+1},p_{k+1})$),
and our ray intersects the segment $AB$.
If the ratio $q_{k+1}/q_k$ is large,
than $OB>>OA$. 

We apply the known procedure of
"noses stretch".
Namely, we define a sequence of points
$B_0=B,\,\overline{OB_i}=\overline{OB_{i-1}}+\overline{OA}$.
One of segments $B_{i-1}B_i$ intersects our ray. 
The triangle  $OB_{i-1}B_i$ is clearly basic,
and, if $a_k$ and $k$ are quite large,
it is $\e$-suitable.
\smallskip

Let domains $\O_1,\,\O_2$ be homothetic triangles
$\O_1=OCD,\,\O_2=\l \O_1\,(\l>0)$, where the line
$CD$ is vertical and points $C,\,D$ lie in the domain
$0<y<x$ ($C$ lower than $D$).

Using Vitali theorem and our lemma,
we may find $\e$-suitable triangles $OA_iB_i\,(i=1,\,2,\,\dots,n)$
such that rays $OC,\,OA_1,\,OB_1,\,OA_2,\,OB_2,\,\dots,OA_n,\,OB_n,\,OD$
go counterclockwise in given order and

$\sum_{i=1}^n \angle A_iOB_i>\angle COD-\e$.

Let $\Delta_{i}$ be the largest triangle,
homothetic (with centre $O$)
to $OA_iB_i$, which is contained in $\O_1\,(i=1,\,2,\,\dots,\,n)$.

Then we may sum up the estimates of type (a3) 
for triangles $\Delta_i,\,\l \Delta_{i}$,
and get the following lower bound:
$$
M(\O_1,\O_2;\,n)\ge 
\left(\l \sum_{i=1}^n 2S(\Delta_i)\right)n^2+o(n^2)
$$
For small $\e$, we have
$S(\Delta_i)\ge c(\e) S(\angle A_iOB_i\cap \triangle OCD)$,
where $c(\e)\rightarrow 1$ if $\e\rightarrow 0$.

So, we let $\e$ tend to zero and for given domains
$\O_1,\,\O_2$ we get a lower bound of type (a3):
$$
M(\O_1,\,\O_2;\,n)\ge \zeta(2)^{-1}\int_0^{\pi} 
l(\O_1,\f)\cdot l(\O_2,\f) d\f\    n^2+o(n^2)
\eqno(a4)
$$

Consider points $C_1$ and $D_1$, 
in which line $CD$ meets X-axis and line
$x=y$, respectively. Consider triangles
$\triangle_1=OC_1C,\,\triangle_2=OCD,\,
\triangle_3=ODD_1,\,\triangle_0=OC_1D_1$ 
and triangles
$\l \triangle_i$.
We have
$$
M(\triangle_0,\l \triangle_0;\,n)\le 
\sum_{i=1}^3\sum_{j=1}^3 M(\triangle_i,\l \triangle_j;\, n).
$$
LHS has an asymptotics of type (a3), three summands
in RHS (with $i=j$) satisfy the lower bound
(a4). Combining these two observations,
we see that the lower bound is an upper
bound aswell in all three cases, and crossing terms
give a contribution $o(n^2)$.

So, the asymptotics (a3) is gotten
for domains of described type.

Now all "quite good" (for example,
with piecewise-smooth boundary)
domains may be approximated from both sides by sums
and differences of such domains.

Let's now consider slightly more
general problem.
Namely, replace the condition on a pair
of vectors $\x_1\times \x_2=1$ to the condition
$\x_1\times \x_2=m=\const\ne 0$.
Again we start from the same special case, 
when
$$\O_i=\{(x,y):0<y<x\le a_i\}\,(i=1,2,\,a_i>0).$$
It's easy problem to find answer in this case.
Indeed, let's fix the greatest common divisor
GCD$(x_1,x_2)=d|m$ ($x_1,\,x_2$ are abscisses
of vectors $\x_1,\,\x_2$). 
We get
$$d^{-1}\zeta(2)^{-1}n^2a_1a_2+o(n^2)$$
our pairs. Sum up by all divisors of $m$
and get asymptotics
$$
M(\O_1,\,\O_2;\,m;\,n)=\sigma(m)\cdot|m|^{-1}\cdot\zeta(2)^{-1}\int_0^{\pi} 
l(\O_1,\f)\cdot l(\O_2,\f) d\f\  n^2+o(n^2)
\eqno(a3')
$$
The generalization for generic
domains $\Omega_1,\,\Omega_2$ does not differ
from the one for case $m=1$.
\smallskip
Lemma 2 easily follows from the proven
fact.

{\bf Proof of lemma 2.}
Without loss of generality, we assume that
$\O$ is a rectangle $\{0<x<A,\,0<y<B\}$.
In this case, the number of blue points
in a rectangle $N\O$ is a number of pairs
of vectors $(\x,\,y): \x\times \y=m$ such that
$\x\in N\O_x,\,\y\in N\O_y$,
where $\O_x$ and $\O_y$ are domains defined
as 
$$\O_x=\{\x\in \An:[\x]\le A\},\,\O_y=\{\y\in \An: [\x]\le B\}.$$
Applying (3'), we get the desired result.

   \smallskip

{\centerline{\bf References.}}
\smallskip

\item{[J]}\ {\it V. Jarn\'\i k}.
{\rm \" Uber} die Gitterpunkte auf konvexen Kurven. 
Math. Z. 1926. Bd.~24. S.~500--518                                        

\item{[P]}\ {\it A. Plagne}. A uniform version of {\rm Jarn\'\i k's} theorem.
Acta Arith., {\bf 57}, no.3 (1999), 255-267.
 
\item{[EM]}  {\it A. Eskin, C. McMullen}.            
Mixing, counting and equidistribution in Lie groups.
Duke Math. J. 71(1993), 181-209.  

\item{[Ve]}\ {\cit A. Vershik}. 
Limit shape of convex polygons, Funkc. anal i pril.
{\bf 28} (1994), 13--20. (Russian) 

\item{[Ba]}\ {\it I. Barany}.
The Limit Shape of Convex Lattice Polygons. 
Discrete and Computational Geometry {\bf 13} (1995), 279-295.

\item{[BP]}\ {\it E. Bomberi and J. Pila}. The number of integral points
on arcs and ovals, Duke Math. J. 59 (1989), 337-357

\item{[SD]}\ {\it H. P. F. Swinnerton-Dyer}. The number of 
lattice points on a convex curve. J. Number Theory 6 (1974), 128--135. 

\item{[Kh]}\ {\cit A. Ya. Khinchin}. 
Continued fractions. M., Fizmatgiz (1961). (Russian)  

\item{[F]}\ {\cit J. Favar}. 
Kurs lokalnoy differencialnoy
geometrii M., IL, 1961. (Russian)
 
             \end